\documentclass{article}
\usepackage{amsmath}
\font\tengoth=eufm10
\font\sevengoth=eufm7
\font\fivegoth=eufm5
\newfam\gothfam
\textfont\gothfam=\tengoth \scriptfont\gothfam=\sevengoth
\scriptscriptfont\gothfam=\fivegoth

\def\blacksquare{\hbox to .60em{\vrule width .60em height .60em}}

  \font\bb=msbm10 
\catcode`\@=12  

\def\é{\'e}
\def\{\`e}
\def\?{\`a}
\def\{\`u}
\def\{\c c}
\def\hb {\hfil \break}
\def\n {\vskip 0.2cm \noindent }
\def\scirc{\,{\raise 0.8pt\hbox{$\scriptstyle\circ$}}\,}
\def\ins{\,{\raise 0.2cm \hbox{ $\scriptstyle \circ$}}\,}

\def  \é{\'e}
\def\è{\`e}
\def\à{\`a}
\def\ù{\`u}
\def\ç{\c c$\!\!\!$}

\date{}
 
\begin{document}
  
 \centerline {\bf  Des s\éries de tissus ordinaires de rang maximum}
 \centerline {\bf en toute dimension}
 
 \bigskip
 
 \centerline{\bf   J.P. Dufour\footnote{Jean-Paul Dufour, ancien professeur \à l'Universit\é de Montpellier II,  
  1 rue du Portalet, 34820 Teyran, France\hb  email : dufourh@netcourrier.com} et  D. Lehmann\footnote 
     {Daniel Lehmann, ancien professeur \à l'Universit\é de Montpellier II,   4 rue Becagrun,  30980 Saint Dionisy, France\hb  email : lehm.dan@gmail.com}}
     
 
 
 \section{Introduction}

 On pose :
 $ c(n,h) := \begin{pmatrix}n+h-1\\h \end{pmatrix}$
  
 \n
  \hskip 3cm (dimension de l'espace des polyn\^omes homog\ènes de degr\é $h$ en $n$ variables).
   
\n Pour tout 
 entier $n$ $(n\geq 2)$,
et tout entier $d$ $(d\geq n)$, 
  on note :
 
 $k_0$ l'entier tel que $c(n,k_0)\leq d\leq c(n,k_0+1),$
 
  $ \pi'(n,d) := 
  k_0.d-c(n+1,k_0)+1\ ,$\n \hskip 1.8cm $=\sum_{h=1}^{k_0}\bigl(d-c(n,h)\bigr),$
 
 \indent  
 (rang maximum des $d$-tissus ordinaires, au sens de [CL][DL], sur une vari\ét\é de dimension $n$,)
 
 $\rho(n,k_0) := \pi'\bigl(n,c(n,k_0)\bigr).$

 \n {\it Convention :} 
Pour tout couple $(p,q)$ d'entiers, $(p\geq 0)$, on   pose : $\begin{pmatrix} p\\q \end{pmatrix}=0$ si
 $q<0$ ou $q>p$.
 
 Dans cet article, nous d\éfinissons  une  famille infinie 
 $\bigl({\cal W}(n,E)\bigr)_n$ de $c(n,k_0)$-tissus holomorphes de codimension un, en toute dimension $n$,  \à partir d'un certain ensemble fini $E$ de fonctions. Nous montrons qu'il suffit de v\érifier un nombre \emph{fini} de conditions pour que ces tissus soient tous ordinaires, et tous de rang maximum $\rho(n,k_0)$. En outre, ces conditions peuvent \^etre test\ées pratiquement sur ordinateur \à l'aide du programme publi\é dans [DL], avec un temps de calcul qui reste acceptable pourvu que $k_0$ ne soit ``pas trop grand''.

 \subsection{D\éfinition des tissus ${\cal W}(n,E)$}
 
 \n {\bf Lemme 1 :} {\it Quels que soient $n$ et $h$, la formule  suivante  est  v\érifi\ée :
 $$c(n,h)=\sum_{k=1}^{h}\begin{pmatrix}h-1\\k-1 \end{pmatrix}.\begin{pmatrix}n\\k \end{pmatrix}.$$}
 
  \n {\it D\émonstration :}     
  L'espace ${\cal L}_h(n)$ des polyn\^omes homog\ènes  de degr\é $h$ en $n$ variables $(x_1,\cdots,x_n)$ se d\écompose en somme directe des sous-espaces ${\cal L}_h(k,n)$ engendr\és par les mon\^omes     dans lesquels  interviennent   exactement $k$ des $n$ variables. Or la dimension de ${\cal L}_h(k,n)$ est  \égale  \à $\begin{pmatrix}h-1\\k-1 \end{pmatrix}.\begin{pmatrix}n\\k \end{pmatrix}$ puisque  $\begin{pmatrix}h-1\\k-1 \end{pmatrix}$   est aussi la dimension $c(k,h-k)$    de l'espace des polyn\^omes homog\ènes  de degr\é $h$ en $k$ variables $(x_1,\cdots,x_k)$  qui sont multiples du mon\^ome $x_1.x_2\cdots.x_k$. Le lemme en r\ésulte.
   
   \rightline{QED}
   
\pagebreak
 
A partir de maintenant,   l'entier $k_0$  $(\geq 2)$  est fix\é une fois pour toutes. 

\n Soit $T_k$ un tissu sur un ouvert de {\bb C}$^k$. Si $n$ d\ésigne un entier $>k$, et $I$ un  multi-indice\break  $I=(1\leq i_1<i_2<\cdots<i_k\leq n)$, la projection  $\pi_I^n:(x_1,x_2,\cdots,x_n)\mapsto (x_{i_1},x_{i_2},\cdots,x_{i_k})$ de {\bb C}$^n$ sur {\bb C}$^k$ permet de transposer    $T_k$ en un tissu $(\pi_I^n)^*(T_k)$ sur un ouvert de {\bb C}$^n$ :  si $T_k$ est d\éfini localement par la  famille  d'int\égrales premi\ères  $\bigl(u_b(x_1,\cdots,x_k)\bigr)_{ b}$,  $(\pi_I^n)^*(T_k)$ admet    $\bigl(u_b\scirc \pi_I^n\bigr)_b$ comme famille d'int\égrales premi\ères.

 \n Dans la suite,   $T_1$ d\ésignera   le feuilletage par points sur un ouvert de {\bb C} (admettant  comme int\égrale premi\ère la fonction identit\é $x\mapsto x$, ou n'importe quelle autre fonction holomorphe r\éguli\ère de la variable $x$). Si $\pi_\lambda^n:(x_1,x_2,\cdots,x_n)\mapsto x_\lambda$ d\ésigne la $\lambda$-\ème projection, $(\pi_\lambda^n)^*(T_1)$ est donc  le feuilletage $(x_\lambda=cte)$.

   \n {\bf D\éfinition :}
 {\it On appellera  ensemble de tissus \emph{$k_0$-\équilibr\é} la donn\ée, 
  pour tout $k=1,\cdots,k_0$ ,  d'un  $c(k,k_0-k)$-tissu $T_k$,  v\érifiant les conditions suivantes :
 
 - les int\égrales premi\ères $u_b$  permettant de d\éfinir $T_k$ font intervenir explicitement toutes les\hb  \indent variables $(x_1,\cdots,x_k)$ ;
 
 
 - si l'on fait varier   $I$ dans l'ensemble $I_k(n)$   des multi-indices $I=(1\leq i_1<i_2<\cdots<i_k\leq n)$, et $k$ de 1 \à $k_0$, la superposition de tous les  tissus $(\pi_I^n)^*(T_k)$ est un tissu  ${\cal W}(n,E)$ sur {\bb C}$ ^n$.
 
 \n $[$Si $n<k_0$, seuls interviennent les tissus $T_k$ tels que  $k\leq n]$.
   }

  Si $E=(T_k)_k$ est un tel ensemble $k_0$-\équilibr\é,    le tissu ${\cal W}(n,E)$ 
d\éfini ci-dessus est un $c(n,k_0)$-tissu en vertu du lemme 1. Il est donc  "calibr\é", selon la terminologie introduite dans [DL], et l'on peut,    d\ès lors qu'il est ordinaire,  tester s'il est ou non de rang maximum, selon que sa courbure est nulle ou non.

  \subsection{Exemples d\éj\à connus}(les tissus $T_k$ \étant d\éfinis par des int\égrales premi\ères\footnote{On se permettra l'abus de notation consistant \à noter une fonction   $u(x_1,\cdots,x_k)$,   au  lieu de $u$ ou de \hb  $ \indent (x_1,\cdots,x_k)\mapsto u(x_1,\cdots,x_k)$  .}) :
  
  - pour $k_0=2$, les tissus ${\cal W}(n,E)$, s'ils sont de rang maximum, sont d'un int\ér\^et limit\é, puisqu'ils sont  tous lin\éarisables et localement du type $T_1=\{x\}$ et  $T_2=\{x+y\}$ ;

  - pour  $k_0=3$ avec  $T_1 = \{x \}$, $T_2 = \{x+y,x-y\}$, 
 $T_3 =\{ x^2+y^2+z^2\}$, on obtient les tissus : \hb $\{x_i$ $(i\leq n)$, $x_i+x_j, \ x_i-x_j$ $(i<j)$, $(x_i)^2+(x_j)^2+(x_k)^2$ $(i<j<k)\}$, ordinaires et de rang maximum $\rho(n,3)$ ;

  -  pour  $k_0=4$, notons $\beta(x,y,z,t)$ le birapport de quatre points distincts dans la droite projective complexe ;   
  les tissus ${\cal W}(A_{0,n+3})$ de Pereira-Pirio ([Pe]) (tissu de Bol pour $n=2$)   sont tous ordinaires et de rang maximum $\rho(n,4)$ ; ils sont obtenus avec 
  $T_1 = \{\beta(x,0,1,\infty) \}$, \hb $T_2 = \{\beta(x,y,0,\infty),\beta(x,y,1,\infty),\beta(x,y,0,1)\}$,  
 $T_3 =\{\beta(x,y,z,0), \beta(x,y,z,1), \beta(x,y,z,\infty) \}$, et \hb
 $T_4 =\{\beta(x,y,z,t)\}$ ;  
 
 \n [La version affine de ces tissus, obtenue  en fixant l'un des 3 points $0,1$ ou $\infty$ (par exemple en prenant les rapports alg\ébriques $\beta(x,y,z,\infty)$ de trois points parmi les $n+2$ nombres  $(x_1,\cdots,x_n,0,1)$), procure \également des tissus ordinaires et de rang maximum $\rho(n,3)$ ; ils sont alors  lin\éarisables].
  
 - pour  $k_0=4$ \également, on obtient  les tissus ${\cal W}B_n  $   de  [DL]  (tissu d\^u \à Pirio ([Pi]) pour $n=2$)  avec   
  $T_1 = \{x \}$, $T_2 = \{x+y,x-y,xy\}$,  
 $T_3 =\{x+y+z, x^2+y^2+z^2, xyz\}$, et\hb 
 $T_4 =\{x+y+z+t\}$ ou $\{x^2+y^2+z^2+t^2\}$ ou $\{xyzt\}$  ; \hb  ils sont encore tous ordinaires et de rang maximum $\rho(n,4)$.
 
 \pagebreak 
 
 \n {\bf Remarque : }
 
 Tous les exemples ci-dessus sont \emph{quasi-sym\étriques}\footnote{Ce concept   d\épend du choix des coordonn\ées locales. Il n'est donc pas intrins\èque, et n'a   de signification que localement.}  au sens suivant : pour tout $k$ $(\leq k_0)$, et     pour toute permutation   $\sigma$ de $\{1,2,\cdots,k\}$, le tissu $T_k$ est $($globalement$)$ invariant par le diff\éomorphisme 
  $$(x_1,x_2,\cdots,x_k)\mapsto (x_{\sigma(1)},x_{\sigma(2)},\cdots,x_{\sigma(k)}).$$ 
Mais cette  propri\ét\é ne sera pas requise en g\én\éral.  
 
 \section{Caract\ère ordinaire des tissus ${\cal W}(n,E)$}

 Pour tout $d$-tissu $W$ de codimension un en dimension $n$ d\éfini localement par les int\égrales premi\ères  $ \bigl(u_i(x_1,\cdots,x_\lambda ,\cdots,x_n)\bigr)_{1\leq i\leq d}$\ ,  et pour tout multi-indice $L=(\ell_1,\cdots,\ell_\lambda,\cdots,\ell_n)$ form\é d'entiers  $\ell_\lambda\geq 0$,  on pose :
 $$c_L^i(W):=\prod_{\lambda=1}^n \Bigl(\frac{\partial u_i}{\partial x_ \lambda}\Bigr)^{\ell_\lambda}.$$
 Ayant ordonn\é les indices $i$ et les multi-indices $L$,  notons  $P_h(W):=\Bigl(\!\Bigl(c_L^i(W)\Bigr)\!\Bigr)_{i,L}$ la matrice de taille $c(n,h)\times d$ d\éfinie pour $1\leq i\leq d$ et $|L|=h$,
 o\ù  $|L|:=\ell_1+\ell_2+\cdots+\ell_n$ d\ésigne le degr\é  de $L$. Rappelons ([CL],[DL]) que le tissu est dit \emph{ordinaire} (au voisinage d'un point de {\bb C}$^n$) si toutes les matrices $P_h(W)$ sont de rang maximum inf$\bigl(d,c(n,h)\bigr)$ au voisinage de ce  point, et qu'il suffit pour cela qu'elles le soient pour $h\leq k_0$, $k_0$ d\ésignant l'entier tel que $c(n,k_0)\leq d<c(n,k_0+1)$. Posons :
 $$P_{h,n} :=P_{h}\bigl({\cal W}(n,E)\bigr).$$
 Supposons aussi chaque tissu $T_k$ d\éfini localement par des int\égrales premi\ères $u_b(x_1,\cdots,x_k)$, \break  $(1\leq b\leq c(k,k_0-k))$ gr\^ace auxquelles on peut d\éfinir les matrices $P_{h}(T_k)$.
 
 Le th\éor\ème suivant permet de ramener la d\émonstration du caract\ère ordinaire de chacun des tissus de la famille infinie 
 $\bigl({\cal W}(n,E)\bigr)_n$ \à la v\érification  d'un nombre fini de conditions :

\n {\bf Th\éor\ème 1 :}

{\it Les quatre   assertions suivantes sont \équivalentes :

 $(i)$   Les tissus ${\cal W}(n,E)$ sont tous ordinaires, quel que soit $n$, $(n\geq 2)$. 
 
 $(ii)$ Le tissu ${\cal W}(k_0,E)$ est  ordinaire.

$(iii)$ La matrice $P_{k_0,k_0} $ est inversible.

$(iv)$ Les matrices carr\ées $P_{d_k}(T_k)$, de taille $d_k\times d_k$, sont toutes inversibles  quel que soit $k$,\hb  \indent $(1\leq k\leq k_0)$, o\ù l'on a    pos\é $d_k:=c(k,k_0-k)$.}

\n {\it D\émonstration :}  Les implications $(i)\Rightarrow (ii) \Rightarrow (iii) \Rightarrow (iv)$ sont \évidentes. On va maintenant montrer   $(i)\Leftrightarrow (iv)$.
 
Commen\ç ons par  ordonner les $c(n,k_0)$ indices $i$ et les $c(n,h)$ indices $L$ de la matrice $P_{h,n}$. Pour cela, donnons 
 nous, pour tout $k=1,\cdots,k_0$, un ordre arbitraire  $O(k,n)$ sur l'ensemble $I_k(n)$ des  multi-indices $1\leq \lambda_1<\lambda_2<\cdots<\lambda_k\leq n $.  Fixons aussi un ordre arbitraire $O'(k)$ sur l'ensemble $E_k$ des    $c(k,k_0-k)$ int\égrales premi\ères $u_b$. Fixons enfin  un ordre arbitraire $O"(k,h)$ sur l'ensemble ${\cal L}_h(k,n)$ des multi-indices $L=(\ell_1,\cdots,\ell_\lambda,\cdots,\ell_{n})$ de degr\é $|L|=h$ pour lesquels il y a  exactement  $k$ des entiers $\ell_\lambda$ qui sont  non-nuls.
 
 \n A chaque fonction int\égrale premi\ère $u_b(x_{\lambda_1},\cdots,x_{\lambda_k})$ de ${\cal W}(n,E)$,  on attribue  le triplet $(k,a,b)$, o\ù 
 
   $k$   d\ésigne  
 le nombre de variables intervenant dans l'expression de la fonction ($1\leq k\leq k_0$), 
 
  $a$  le num\éro d'ordre   de 
 $(1\leq \lambda_1<\cdots<\lambda_k\leq n)$   dans  ${ I}_k(n)$ pour $O(k,n)$,
 
   et $b$   le num\éro d'ordre   de $u_b$ dans $E_k$ pour $O'(k)$,   $\bigl(1\leq b \leq c(k,k_0-k)\bigr)$.

  \n Au   multi-indice $L=(\ell_1,\cdots,\ell_\lambda,\cdots,\ell_{n})$,  on attribue  le quadruplet $(h,k,a,b)$, o\ù 
  
  $h=|L|$,
  
   $k$  d\ésigne celui des indices tels que   $L\in {\cal L}_h(k,n)$ ($1\leq k\leq h$), 
 
  $a$  le num\éro d'ordre   du multi-indice   
 $(\lambda_1<\cdots<\lambda_j<\cdots<\lambda_k)$    dans  ${ I}_k(n)$ pour  $O(k,n)$,  tel que\hb  \indent $\ell_{\lambda_j}\neq 0$ pour tout $j=1,\cdots,k$,  
 
   et $b$   le num\éro d'ordre   de $L$  dans ${\cal L}_h(k,n)$ pour $O"(k,h)$,  $\bigl(1\leq b\leq c(k,h-k)\bigr) $.
   
   \n \n On range alors  les fonctions (resp. les quadruplets de degr\é $h$)  dans l'ordre alphab\étique des triplets  :\hb  $(k,a,b)<(k',a',b')$ si $k<k'$ ou ($k=k',a<a'$), ou 
 ($k=k',a=a',b<b'$). 
   
 \n  On note $P_{h,n}\bigl((k,a),(k',a')\bigr)$ le bloc dans la matrice $P_{h,n}$ form\é avec les indices colonne de la forme $(k,a,\ast)$ et les indices ligne de la forme $(h,k',a',\ast)$. Ces blocs ne sont form\és que de z\éros si $k<k'$ (ou si $k=k'$ et $a\neq a'$).

  \n Les blocs sous-diagonaux n'ayant que des z\éros, la matrice $P_{h,n}$ est de rang maximum $c(n,h)$ ssi  tous les blocs diagonaux $P_{h,n}\bigl((k,a),(k,a)\bigr)$ sont   de  rang maximum $c(h,k_0-h)$.  Mais 
   chaque  bloc diagonal $P_{h,n}\bigl((k,a),(k,a)\bigr)$ est \égal, au nom pr\ès des variables, \à la matrice $P_{h}(T_k)$, ce qui montre que $(i) $ \équivaut au fait que toute les matrices $P_h(T_k)$ sont de rang maximum, (et   entra\^ine  aussi l'implication $(iii)\Rightarrow (iv)$).
   
    Montrons alors qu'il suffit que $P_{d_k}(T_k)$ soit de rang maximum $d_k$ pour que les 
 matrices  $P_{h}(T_k)$ soient toutes  de rang maximum quel que soit $h$, ($h\leq d_k)$. 
 
 Pour $k=1$, il  n'y a rien  \à d\émontrer.  Pour $i=(k,a,b)$ avec $k\geq 2$, tous les  feuilletages $u_i=cte$ 
 sont transverses 
 aux  feuilletages $x_\lambda =cte$, puisque ceux-ci font partie du tissu  ;  par cons\équent toutes les d\ériv\ées partielles $\bigl(u_b(x_{i_1},\cdots,x_{i_k})\bigr)'_{i_j}$ sont non-nulles.  
  
  Soit $a=(i_1<\cdots <i_k)$. A chaque ligne  $L=(\ell_1,\cdots,\ell_\lambda,\cdots,\ell_{n})$ 
   de $P_h(T_k)$ correspondant au quadruplet $(h,k,a,b)$,   associons l'indice  $\lambda_k(L)$, le dernier des $\lambda_j$ tels que  $\ell_{\lambda_j}\neq 0$, 
     et multiplions tous les termes de la colonne $i=(k,a,b)$ de $P_h(T_k)$ par le nombre   $\Bigl(\bigl(u_b(x_{i_1},\cdots,x_{i_k})\bigr)'_{\lambda_k(L)}\Bigr)^{d_k-h}$ : en r\ép\étant l'op\ération pour chaque colonne de $P_h(T_k)$, on obtient la ligne $L'=(\ell'_1,\cdots,\ell'_\lambda,\cdots,\ell'_{n})$ de  $P_{d_k}(T_k)$ pour laquelle $\ell'_\lambda=\ell_\lambda$ si $\lambda<>\lambda_k$, et $\ell'_{\lambda_k}=\ell_{\lambda_k}+d_k-h$. En r\ép\étant maintenant l'op\ération pour chaque ligne de $P_h(T_k)$, on obtient une sous-matrice de $P_{d_k}(T_k)$, qui est de rang maximum si l'on suppose $P_{d_k}(T_k)$  elle m\^eme de rang maximum, et  qui a m\^eme rang que $P_{h}(T_k)$ puisque les nombres  $\bigl(u_b(x_{i_1},\cdots,x_{i_k})\bigr)'_{\lambda_k}$ ne sont pas nuls. Ceci ach\ève de prouver 
   $(i)\Leftrightarrow (iv)$.

  \rightline{QED}
  
   \section{Tissus ${\cal W}(n,E)$ de rang maximum}
   
   L\à encore, la v\érification d'un nombre fini de conditions va nous permettre de conclure pour la famille infinie de tissus.

\n {\bf Th\éor\ème 2 :}

{\it Si tous les tissus ${\cal W}(n,E)$ sont ordinaires  $(n\geq 2)$, il suffit - pour qu'ils soient tous  de rang maximum $\rho(n,k_0)$ quel que soit $n$ - qu'ils le soient 
 pour $ n\leq k_0$.
 
 En outre, leurs relations ab\éliennes peuvent s'exprimer avec au plus $k_0$ variables ; et si $h$ d\ésigne un entier compris entre 2 et $k_0$, la dimension du sous-espace  des relations ab\éliennes   faisant intervenir exactement $h$ des variables $(x_1,\cdots,x_n)$ est   \égale \à $N(h,k_0).\begin{pmatrix}n\\h \end{pmatrix}$ d\ès que $n\geq h$, $N(h,k_0)$ d\ésignant  la dimension du sous-espace  des relations ab\éliennes de  ${\cal W}(h,E)$ faisant intervenir exactement toutes  les $h$ variables $(x_1,\cdots,x_h)$. }

\n {\it D\émonstration :} 

Puisque   les tissus ${\cal W}(n,E)$  sont   ordinaires,   le nombre $\rho(n,k_0)$
 est    une borne sup\érieure de leur  rang,    d'apr\ès [CL] ou [DL]. Il nous suffit donc de montrer que ce nombre  est aussi une borne inf\érieure quel que soit $n>k_0$, d\ès lors que c'est vrai pour $n\leq k_0$.

D\éfinissons pour cela l'entier $N(n,k_0)$ par r\écurrence sur $n$ ($n\geq 2$) en posant :
 
 $N(2,k_0):= \rho(2,k_0),$
 
 $N(n,k_0):=\rho(n,k_0)-\sum_{h=2}^{n-1}N(h,k_0).\begin{pmatrix}n\\h \end{pmatrix}$ pour $n>2$. 
 
Puisque ${\cal W}(2,E)$ est de rang maximum $\rho(2,k_0)$,   ${\cal W}(n,E)$  admet  exactement $\rho(2,k_0).\begin{pmatrix}n\\2 \end{pmatrix}$ relations ab\éliennes  ind\épendantes faisant intervenir explicitement deux des $n$   coordonn\ées    $(x_1,x_2,\cdots,x_n)$. En particulier, 
 ${\cal W}(3,E)$, puisqu'il est de rang maximum $\rho(3,k_0)$,   admet exactement  $\rho(3,k_0)-\rho(2,k_0).\begin{pmatrix}3\\2 \end{pmatrix}$    relations ab\éliennes  ind\épendantes (soit $N(3,k_0)$)  faisant intervenir explicitement les  trois    coordonn\ées; et   ${\cal W}(n,E)$  admet  exactement   $N(3,k_0).\begin{pmatrix}n\\3 \end{pmatrix}$ relations ab\éliennes  ind\épendantes faisant intervenir explicitement trois des $n$   coordonn\ées    $(x_1,x_2,\cdots,x_n)$ d\ès que $n\geq 3$. De proche en proche, on montre que pour tout $h$ ($2\leq h\leq k_0$), ${\cal W}(n,E)$  admet  exactement $N(h,k_0).\begin{pmatrix}n\\h \end{pmatrix}$ relations ab\éliennes  ind\épendantes faisant intervenir explicitement $h$ des $n$   coordonn\ées    $(x_1,x_2,\cdots,x_n)$, d\ès que $n\geq h$.

  Ainsi, pour $n\geq k_0$,  ${\cal W}(n,E)$ admet   au moins  $\sum_{h=2}^{k_0}N(h,k_0).\begin{pmatrix}n\\h \end{pmatrix}$  relations ab\éliennes ind\épen-\hb dantes. Or ce nombre est pr\écis\ément \égal \à $ \rho(n,k_0)$, comme l'\énonce le lemme ci-dessous. 
  
  \rightline{QED}
   
  \n {\bf Lemme 2 :} {\it La formule  suivante  est  v\érifi\ée  pour $n\geq k_0 $ :
 $$ \rho(n,k_0)=\sum_{h=2}^{k_0}N(h,k_0).\begin{pmatrix}n\\h \end{pmatrix}$$
 ou, de fa\ç on \équivalente, $N(h,k_0)=0$ pour $h>k_0$.}
 
 \n {\bf Remarque :} La formule est en fait v\érifi\ée pour tout $n\geq 2$ ; mais  pour $n\leq k_0$, ce n'est que la d\éfinition de $N(n,k_0)$.
 
 \n {\it D\émonstration du lemme :} Les deux termes de l'\égalit\é sont des polyn\^omes de degr\é $k_0$ par rapport \à la variable  $n$,  
  qui prennent les m\^emes valeurs en $k_0+1$ points distincts (d'abord pour tout entier $n$ compris entre 2 et $k_0$ d'apr\ès la remarque ci-dessus, et aussi pour $n=0$ et $n=1$ o\ù  ils  s'annulent): ils sont donc  n\éc\éssairement \égaux. 
  
   \rightline{QED}
   
  \section{Nouveaux exemples}

   - pour $k_0=3$, les   tissus ${\cal W}(n,E)$ construits avec les ensembles $E$ suivants 
   
   \indent $T_1=\{x \}$, $T_2=\{x+y, xy \}$, et $T_3=\{x+y+z \}$ ou  $\{xyz \}$,
   
    \indent $T_1=\{x \}$, $T_2=\{ xy, \frac{(x-1)(y-1)}{(x+1)(y+1)} \}$, et $T_3=\{x+y+z \}$ ou  $\{\frac{(x-1)(y-1)(z-1)}{(x+1)(y+1)(z+1)} \}$,
    
     \indent $T_1=\{x \}$, $T_2=\{ x+y, \frac{1}{x} +\frac{1}{y}\}$, et $T_3=\{x+y+z \}$ ou  $\{\frac{1}{x} +\frac{1}{y}+\frac{1}{z} \}$ 
     
     \n sont ordinaires et de rang maximum  : il suffit de le v\érifier \à l'ordinateur avec le programme de [DL] pour $n=2$ ou $3$. [Les tissus obtenus pour $n=2$ sont en fait, en plus des tissus lin\éaires,  tous\footnote{C'est aussi avec le programme de [DL] qu'on peut v\érifier le fait qu'on les a tous.} les\hb  4-tissus planaires \emph{presqu'hexagonaux}, c'est-\à-dire isomorphes \à des  tissus de rang maximum du type  $\{x, y, f_1(x,y), f_2(x,y) \}$   tels que les 3-sous-tissus $\{x, y, f_1(x,y) \}$ et $\{x, y, f_2(x,y) \}$ soient tous deux hexagonaux]. 
     
    \n  - pour $k_0=4$, les   tissus ${\cal W}(n,E)$ construits avec l' ensemble  $E$ suivant 
     
    \indent $T_1=\{x \}$, $T_2=\{x+y, x-y, exp(x)+exp(y) \}$,  $T_3=\{x+y+z, x+y-z, exp(x)+exp(y)+exp(z) \}$,\hb \indent et  $T_4=\{x+y+z+t \}$

       \n sont ordinaires et de rang maximum  : il suffit de le v\érifier \à l'ordinateur avec le programme de [DL] pour $n=2,3$ ou $4$. [ Pour $n=2$,  on obtient  en fait l'un des tissus  de la liste des 5-tissus planaires quasi-lin\éaires de rang maximum de Pirio-Tr\épreau ([PT])].

  Donnons nous maintenant une fonction holomorphe r\éguli\ère $f(x_1,\cdots,x_{k_0})$ de $k_0$ variables sur un ouvert de {\bb C}$^{k_0}$ (ou plus g\én\éralement de ({\bb P}$_1)^{k_0}$) , et un ensemble 
  $\{m_1,\cdots,m_{k_0-1}\}$ de $k_0-1$ points distincts dans {\bb C}  ou  {\bb P}$_1$. 
   Pour tout $k$ ($1\leq k\leq k_0$), d\éfinissons $T_k$ comme l'ensemble des $c(k, k_0-k)$ fonctions 
  $$f(x_1,\cdots,x_{k},m_{i_1},m_{i_2},\cdots,m_{i_{k_0-k}})$$ obtenu en faisant varier le multi-indice 
  $ (1\leq i_1< i_2<\cdots<i_{k_0-k}\leq k_0-1)$ de toutes les fa\ç ons possibles. Supposons l'ensemble   $E=(T_k)_k$ ainsi d\éfini $k_0$ \équilibr\é.  
 
 \n Exemples de tels tissus qui soient ordinaires et de rang maximum : 
 
 - avec $k_0=4$, $f(x,y,z,t)=\frac{(x-z)(y-t)}{(y-z)(x-t)} $, et $\{m_1,m_2,m_{3}\}=\{0,1,\infty\}$, on obtient les   tissus ${\cal W}(A_{0,n+3})$ de Pereira-Pirio.
 
 - avec $k_0=3$, $f(x,y,z)=\frac{x-z }{y-z } $, et $\{m_1,m_2 \}=\{0,1\}$, on obtient les tissus lin\éarisables, constituant la version affine des pr\éc\édents. 
 
  \n {\it Probl\ème :}   
 Trouver plus g\én\éralement des  conditions sur $f$ et sur $\{m_1,\cdots,m_{k_0-1}\}$ pour que la s\érie des tissus ${\cal W}(n,E)$ ainsi d\éfinie
 v\érifie les conditions des th\éor\èmes 1 et 2.

  \n {\bf   R\éf\érences   : }
  

  
\noindent [CL] V. Cavalier et D. Lehmann, Ordinary holomorphic webs of codimension one,  arXiv 0703596 v2[mathDS],  2007, et  Ann. Sc. Norm. Super.  Pisa, cl. Sci (5), vol XI (2012), 197-214.

\noindent [DL] J.P. Dufour et D. Lehmann, Calcul explicite de la courbure des tissus calibr\és ordinaires,  arXiv 1408.3909v1[mathDG], 18/08/2014.







\noindent [Pe]   J. V. Pereira, Resonance webs of  hyperplane arrangements, Advanced studies in Pure Mathematics  99, 2010, 1-30.

\noindent [Pi]  L. Pirio, Sur les tissus planaires de rang maximal et le probl\ème de Chern, note aux CRAS, s\ér. I, 339 (2004), 131-136.

\noindent [PT]  L. Pirio et J.M. Tr\épreau, Abelian functional equations, planar web geometry and polylogarithms, Selecta Mathematica, N.S., 11, n° 3-4, 2005,  453-489.

  \end{document}